\pdfoutput=1
\documentclass{amsart}
\usepackage{amsmath}
\usepackage{amsfonts}
\usepackage{amssymb}
\usepackage{epsfig}
\usepackage{amscd}
\usepackage[all]{xy}
\usepackage{color}

\theoremstyle{plain}
\newtheorem{theorem}{Theorem}

\theoremstyle{definition}

\newcommand\dist{\operatorname{distance}}

\begin{document}

\begin{center}

{\bf Rigorous computations with an approximate Dirichlet domain}

\medskip
{\bf Maria Trnkov\'a} 

{\it Department of Mathematics, University of California, Davis, CA 95616, USA.}
{\it mtrnkova@math.ucdavis.edu}

\medskip

{\bf Abstract} 

\end{center}

In this paper we address some problems concerning an approximate Dirichlet domain. We show that under some assumptions the approximate Dirichlet domain can work equally well as an exact Dirichlet domain. In particular, we consider a problem of tiling a hyperbolic ball with copies of the Dirichlet domain. This problem arises in the construction of the length spectrum algorithm which is implemented by the computer program SnapPea. Our result explains the empirical fact that the program works surprisingly well despite it does not use exact data. Also we demonstrate an improvement in the algorithm for rigorous construction of the length spectrum of a hyperbolic 3-manifold.


\section{Introduction} 

Every hyperbolic 3-manifold $M$ of finite volume can be viewed as a convex finite sided polyhedron in $\mathbb{H}^3$ with pairs of faces identified. This polyhedron is called a {\it Dirichlet domain} and there are many reasons why it is important. It is used in the constructions of the length and ortholength spectra. It is also necessary for an isometry test by using the drilling of the geodesics. The length spectrum is then also used to study the symmetry group of $M$ and Dehn parental test. In practice it is very hard to construct an exact Dirichlet domain but we show that for some problems it is enough to know only an approximate Dirichlet domain. 

Poincare's Fundamental Polyhedron theorem verifies the discreteness of a group of isometries in a hyperbolic space by its action on a polyhedron which satisfies special conditions \cite{Poincare1}, \cite{Poincare2}, \cite{Maskit}. 
The works of \cite{Riley},  \cite{Lip} and \cite{Siegel} describe a rigorous algorithm for the discreteness of a subgroup of hyperbolic isometries of $PSL(2,\mathbb C)$ by verifying the conditions of the Poincare Theorem.  Given a subgroup $G$ as a collection of generating matrices as elements of $PSL(2,\mathbb C)$ they constructed a provably correct Dirichlet domain and showed that under special restrictions on the subgroup $G$ this domain satisfies the conditions of the Poincare's theorem. This implies that the subgroup $G$ is discrete and therefore, the corresponding manifold $\mathbb{H}^3/G$ is hyperbolic.
 
In many cases it is known that a subgroup $\Gamma \subset PSL(2,\mathbb C)$ is discrete, ie. $\Gamma$ is a fundamental group $\pi_1(M)$ of a hyperbolic manifold. For example, this can be proved rigorously using verified computations HIKMOT \cite{HIKMOT}. However, the Dirichlet domain is not known but under some assumption we can still get some important knowledge about it.

Everywhere in this paper if not mentioned explicitly we assume that $M$ is a finite volume hyperbolic 3-manifold with a fundamental group presentation $\Gamma$ and a Dirichlet domain $D$ centered at a point $x$ which corresponds to $\Gamma$. Let $D'$ be a finite sided polyhedron obtained by constructing a Dirichlet domain centered at $x$ for $\Gamma$ but in a finite number of steps such that all generators and their inverses in $\Gamma$ correspond to some of its face-pairing relations. We will call $D'$ an {\it approximate Dirichlet domain}. Notation $D_x$ will be used to specify a Dirichlet domain centered at the basepoint $x$. Otherwise we use $D$. 

Our ultimate goal is to show that if a subgroup of hyperbolic isometries $\Gamma\subset PSL(2,\mathbb C)$ is discrete then an associated exact Dirichlet domain $D$ can be replaced by an approximate Dirichlet domain $D'$. Some properties of $D$ and $D'$ are closely related if a difference of their volumes is small enough and a spine radius of $D'$ is close to a minimal. For example, it applies to the case of tiling a ball in a hyperbolic space an computation of an injectivity radius. A main theorem of this paper tells that $D$ can be replaced by $D'$ and it will not affect the rest of the construction of the length spectrum algorithm.  

\medskip

\begin{theorem} 
If a difference $\Delta V$ of volumes  $D'$ and $D$ is less than some small positive $\epsilon$ and the basepoint $x$ minimizes a spine radius $r'$ of $D'$
then every translation of $D$ inside a ball $B(x,R)\subset \mathbb H^3$  can be obtained by translations of $D'$  without leaving the ball $B$ by using an algorithm of adding neighbors to its faces starting from the domain $D'_x$.
\end{theorem}

$D'$ has a finite volume and is a superset of $D$, not necessarily a proper superset.  Therefore it is a Siegel set \cite{Borel}. 
We use a word {\it covering} instead of {\it tiling} for $D'$ because its covering will have overlaps in $\mathbb H^3$ which are not allowed for tiling. Cayley graphs $C$ and $C'$ are associated with domains $D$ and $D'$. Their vertices are elements of $\Gamma$. Their generating sets are elements of $\Gamma$ which correspond to face-pairings of $D$ and $D'$ respectively. Such a Cayley graph is a lift of a one-skeleton dual to a domain. The Theorem 1 can be formulated in terms of Cayley graphs.

\medskip

 \noindent {\bf Corollary:} {\it Cayley graph $C'$ is connected inside a ball $B(x,R)$ under conditions of Theorem 1.}
  
 \medskip

The paper is organized as follow:  Some required background is presented in Section 2. Section 3  starts with the motivation on the study of the length spetrum. Then it describes algorithms of the construction of the Dirichlet domain, tiling a ball by its copies and a general idea of a length spectrum algorithm. For a fixed basepoint and a given presentation of a fundamental group there is a natural way how to build a Dirichlet domain by taking an intersection of all half spaces. The open question here is how to get all necessary half spaces. 
 
A proof of the main theorem uses Lemmas 1 and 2  and it is described in Section 4. First, we show that the Cayley graph $C'$ is connected in the universal cover and it has the same set of vertices as $C$. Then we show that under some restrictions on $\Delta V$ and the spine radius $r'$, $D'$ cannot contain too many extra faces which do not appear in $D$ and it allows to cover the ball $B(x,R)$ by the same tiling algorithm as for the exact Dirichlet domain. Section 5 describes how to compare rigorously if two elements of a fundamental group are the same or not. This can be applied for further implementations of a rigorous algorithm for length spectra.

\medskip

{\bf Acknowlegments:} The author thanks Robert Haraway, Neil Hoffman, Matthias Goerner, Joel Hass and Misha Kapovich for many helpful conversations  and useful comments. The author also thanks the IAS and Monash University where part of this work has been done. The author acknowledges support from U.S. National Science Foundation grants DMS 1107452, 1107263, 1107367 "RNMS: GEometric structures And Representation varieties" (the GEAR Network) and The Robert Bartnik Visiting Fellowship from Monash University.

\medskip

\section{Background} 

We start this section with reviewing some basic definitions and facts from hyperbolic geometry and topology of 3-manifolds that will be used later.

\medskip

{\bf Definition}: A hyperbolic 3-manifold $M$ is a quotient $\mathbb H^3/\Gamma$ of a three-dimensional hyperbolic space $\mathbb H^3$ by a subgroup $\Gamma$ of hyperbolic isometries $PSL(2,\mathbb C)$ acting freely and properly discontinuously.

\medskip

\noindent The subgroup $\Gamma$ is isomorphic to the fundamental group $\pi_1(M)$.

\begin{theorem} [Mostow-Prasad Rigidity, \cite{M},\cite{Th}]
If $M_1$ and $M_2$ are complete finite volume hyperbolic $n$-manifolds, $n>2$, any isomorphism of fundamental groups $\varphi: \pi_{1}(M_{1}) \rightarrow \pi_{1}(M_{2})$ is realized by a unique isometry.
\end{theorem}

The meaning of this theorem is very important. It tells that geometric invariants of manifolds (e.g. volumes, geodesic lengths) are topological invariants. 
Every element $\gamma \in\Gamma$ of the fundamental group corresponds to a closed geodesic $g\subset M$ of a manifold. On the other hand every pre-image of a geodesic $g\subset M$ in the universal cover $\mathbb H^3$ is preserved by some element of the fundamental group $\gamma$ or its conjugates $q^{-1}\,\gamma \,q$, where $q$ is any element of the fundamental group $\Gamma$. We often identify a geodesic with its conjugacy class in the fundamental group. This relation allows to naturally define a length of a closed geodesic and we use the notation $l(\gamma)$ equivalently to $l(g)$. 

\medskip

 {\bf Definition}: \emph{Complex length} $l(\gamma)$ of a closed geodesic $g$ in a hyperbolic 3-manifold is a number $\lambda+i\theta$, where $\lambda$ is a minimal distance of the transformation $\gamma$, and  $\theta$ is the angle of rotation incurred by traveling once around the axis of $\gamma$, modulo $2\pi$.

\smallskip

 {\bf Definition}: \emph{Length spectrum} $L(M)$ of a hyperbolic 3-manifold $M$ is a set of all complex lengths of closed geodesics in $M$ taken with multiplicities:
 $$
 L(M)=\{l(\gamma)\subset \mathbb{C}|\forall\gamma\in\Gamma\}.
 $$
This is an infinite discrete and ordered set.

\medskip

Another way of thinking about a manifold is via a Dirichlet domain.

\medskip

{\bf Definition:} \emph{Dirichlet domain} of $\Gamma$ with a basepoint $x\in \mathbb H^3$ is the subset  
$$
D_x(\Gamma)=\{y\in X\, |\, d(x,y) \leq d(g(x),y), \forall g\neq 1\in\Gamma\}.
$$
For finite volume manifolds this is a finite sided polyhedron with pairs of faces identified \cite{BP},  \cite{Bo}.
 

\medskip

\section{Length spectrum} 

This section highlights the importance of length spectra of hyperbolic 3-manifolds and their relation to an approximate Dirichlet domain construction. We review some relevant definitions and the algorithm for the construction of the length spectrum. This algorithm was implemented in the computer programs SnapPea, Snap, SnapPy \cite{SnapPy} and {\it Mathematica} package Ortholength.nb.

The length spectrum is a topological but not a complete invariant for hyperbolic 3-manifolds. The first examples were introduced by M.-F.\,Vigneras \cite{Vig} who showed how to construct iso-spectral but non-isometric, commensurable hyperbolic 3-manifolds. There is still an open question if there are iso-spectral manifolds which are not commensurable, see \cite{Reid}, \cite{Futer}. R.~Meyerhoff \cite{Meyerhoff} proved that finite subsets of length, ortholength and basin spectra define a closed hyperbolic 3-manifold uniquely up to an isometry.

Constructions of the ortholength and the basin spectra depend on a length spectrum. Note also that the basin spectrum is a structured subset of the ortholength spectrum. The length spectrum can be used to compute a subgroup of the group of isometries of a closed hyperbolic 3-manifold $M$ \cite{Kojima}. It requires to drill out a complete simple geodesic link in $M$. This drilling results in a complete hyperbolic manifold of finite volume \cite{Sakai} and it can be used for testing if two closed manifolds are isometric \cite{HW}.

Dehn Parental Test asks if $N$ is a Dehn filling of $M$ where $M$, $N$ are orientable 3-manifolds which admit complete hyperbolic metrics of finite volume on their interiors. The first practical method for determining Dehn filling heritage was described by C.D.Hodgson, S.P.Kerckhoff \cite{HoKe}, and more recently implemented by R.Haraway \cite{Ha}. Dehn parental test for hyperbolic 3-manifolds reduces to rigorous calculations of length spectra and also volumes of manifolds, cusp area, slope length and an isometry test.

Computer program SnapPea \cite{JW} written in the C language calculates some topological and numerical invariants of a hyperbolic 3-manifold $M$ including the length spectrum. A manifold $M$ can be given by its ideal triangulation, as a knot or a link complement in $\mathbb S^3$ or a Dehn filling of a cusped manifold. It finds a fundamental group of a manifold first, constructs a Dirichlet domain using the fundamental group and then tiles a ball in $\mathbb H^3$ by copies of the Dirichlet domain. The details of this construction can be found in \cite{HW}, SnapPea's source code. The algorithm is theoretically rigorous but in practice it works with real numbers and cannot avoid round-off errors.

Information about length and ortholength spectra played a crucial role in a proof of the Exceptional manifolds conjecture \cite{GT}. For that purpose authors have implemented the algorithm for a computation of the length spectrum of an oriented hyperbolic 3-manifold in {\it Mathematica} called {\it Ortholength.nb} \cite{GT_code}. The package takes data from Snap \cite{Snap} - a computer program that studies arithmetic and numeric invariants of hyperbolic 3-manifolds and is based on a computer program SnapPea and a number theory package Pari. Ortholength.nb computes the length and ortholength spectra with round-off errors and therefore, it gives a rigorous answer if an input is correct and is given with a big enough precision. In practice, the given precision is over 60 digits. The authors compared results obtained from SnapPea and Ortholength.nb for tens of manifolds with real length cutoff-up to 7 with perfect agreement. The first difference in the spectra appeared when a precision was reduced to less than 10 digits.

We give a sketch of the algorithm for the Dirichlet domain construction as in SnapPea. A point $x\in\mathbb H^3$ is chosen as a basepoint.  $D'$ is obtained as an intersection of half-spaces $H_g$ for all words $g\in\Gamma$ of length $n\in \mathbb N$. If no new face was built for all words of length $n+1$ then it stops and checks vertices, volume, faces, Euler characteristic of $D'$. It does not compute a round-off error but it makes all computations with a predefined error $\epsilon$. It can miss some very narrow faces far from $x$. The SnapPea's construction of a Dirichlet domain is very neat but there is no guarantee that an exact Dirichlet domain was obtained because of the predefined error $\epsilon$. Therefore, it may finish with a construction of an approximate domain. After an approximate Dirichlet domain was constructed it computes a {\it spine radius} of $D'$ - a maximum of all minimal distances from $x$ to an edge of $D'$ over all edges of $D'$. Then the algorithm minimizes the spine radius by slightly moving the basepoint $x$ in its neighborhood. This process is repeated until a minimal spine radius is obtained up to an error $\epsilon$. A result is a more symmetric Dirichlet domain.

The length spectrum algorithm is best described in \cite{HW}. Each transformation $g\in\Gamma$ corresponds to a translation along a fixed geodesic. A translation of $D$ is a $gD$, where $g$ is an isometry of covering transformation $\Gamma$. All geodesics of real length $\lambda$ correspond to group elements $g\in\Gamma$ which can be constructed from products of face pairing relations of Dirichlet domain. All of them move $x$ by a distance less than $R$ which depends on a cut-off real length $\Lambda$ and a spine radius $r$:
$$  
R=2\cosh^{-1}\left(\cosh(r)\,\cosh\frac{\lambda}{2}\right).
$$
Tiling of a ball $B(x,R)$ will produce a {\it big list} of geodesics which might contain geodesics of real length zero or bigger than $\Lambda$, not intersecting a spine, duplicates (powers or conjugates). A {\it small list} of geodesics with correct multiplicities will be obtained after removing all duplicates and geodesics which fail to satisfy some of the constraints.

\section{Rigorous tiling by an approximate Dirichlet domain} 

In this section we show that even if a construction of a Dirichlet domain is not an exact one we can still get an exact tiling of a ball in the universal cover with overlappings. It means the big list of geodesics from Section 3 will not miss any element $g\in\Gamma$ of real length less than $\lambda$.

Fix some point $x$ in the universal cover to be a basepoint of a Dirichlet domain. Then a presentation of a fundamental group $\Gamma$ defines a unique Dirichlet domain D. For the construction of the algorithm we need to tile a ball $B(x,R)$ of some radius $R$ centered at $x$ by copies of D. In practice it is very hard to construct an exact Dirichlet domain \cite{Lip}, \cite{Siegel}. One needs to define all vertices with round-off error or with interval arithmetic. 

We start to construct a Dirichlet domain $D'$ for a fundamental group $\Gamma$ as an intersection of half-planes $H_g$ for finitely many $g\in\Gamma$ including generators and their inverses. We do not know in advance how many steps are involved and if we obtained all vertices, edges and faces of the Dirichlet domain D. Thus the domain $D'$ is a superset of the domain $D$. Faces of $D'$ will contain some faces of $D$ but $D$ may have more faces which lie inside $D'$. If some generator of $\Gamma$ does not correspond to a face of $D'$ then we can change a presentation of $\Gamma$ by replacing a removed generator with a new generator. The new generator maybe the one  which removed the old one. Thus we can always construct $D'$ in such a way that it always contains faces corresponding to generators and their inverses of $\Gamma$.  

\medskip

{\bf Lemma 1:} {\it The injectivity radius $\rho'$ of $D'$ is rigorous if a difference $\Delta V$ between volumes $D$ and $D'$ is small enough.}

\medskip

{\bf Proof:}  By construction an injectivity radius $\rho'$ of an approximate Dirichlet domain $D'$ is always greater or equal to the injectivity radius of the Dirichlet domain $D$. Assume that a face of $D'$ at distance $\rho'$ to the basepoint $x$ is removed and there is a new edge of $D$ which is closer to the basepoint $x$ than $\rho'$. Then two images $g(x)$ and $f(x)$ of the basepoint $x$ corresponding to these two faces: the old and the new one, should belong to different translations of $D$ and a distance between them should be more or equal to $2\rho$ where $\rho$ is a new injectivity radius obtained from constructing the new face. See Figure 1. Then it will give a volume difference between $D$ and $D'$ more than $\Delta V$, where $\Delta V$ can be small enough. $\square$

\begin{figure}
  \includegraphics[scale=0.5]{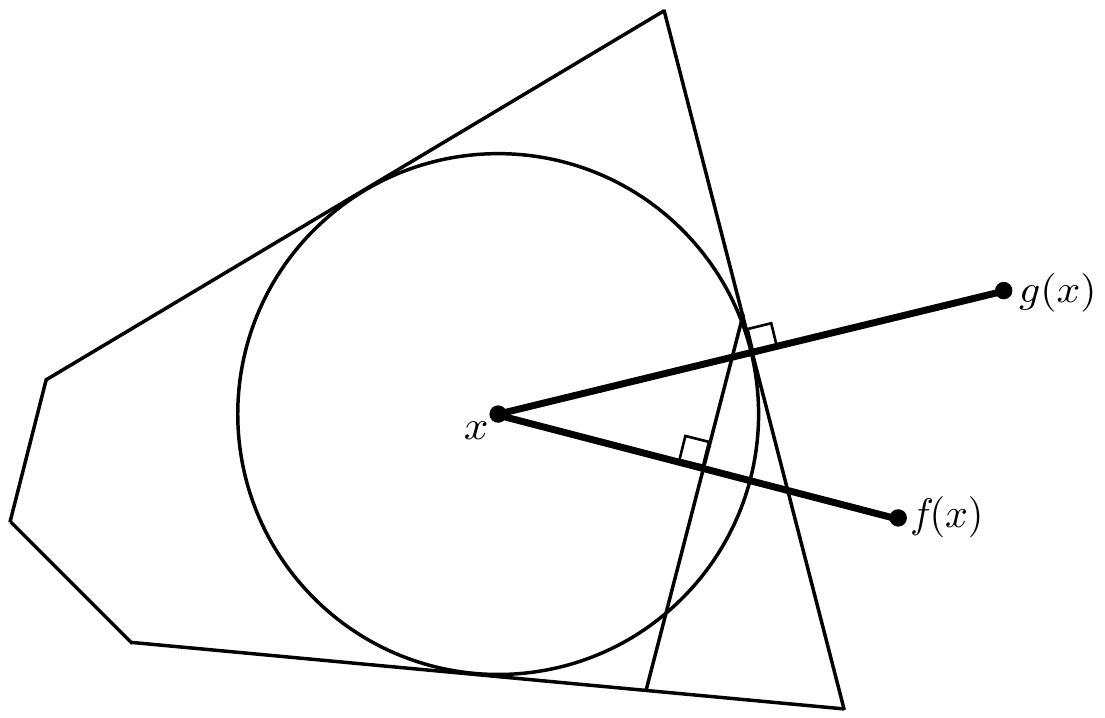}\\
  \caption{Distance between images of the basepoint $x$.}
\end{figure}


If needed one can show in a similar way that a spine radius of $D$ is bigger or equal to a spine radius of $D'$ if difference of volumes $D$ and $D'$ is small enough.

We know that $D$ is a Dirichlet domain of a hyperbolic 3-manifold and tiles the whole universal cover without gaps and overlappings. We replace $D$ by $D'$ and start tiling $\mathbb H^3$ by adding neighborhoods of $D'$ along its faces. It will not be a tiling but a covering of $\mathbb H^3$ because we get all translations of $D'$ without gaps but with overlappings if $D$ is a proper subset of $D'$. This tiling by neighbors can be considered as moving on a Cayley graph along edges. Remind that all generators and their inverses of $\Gamma$ correspond to some faces of $D'$. 

\medskip

{\bf Lemma 2:} {\it Cayley graphs $C$ and $C'$ have the same sets of vertices. Their sets of edges are different if $D$ is a proper subset of $D'$.}

\medskip

{\bf Proof:} Every translation of $D\subset\mathbb H^3$ is obtained by some word $g$ consisting from generators and their inverses of $\Gamma$. So, considering all possible words of $\Gamma$  will map $D'$ onto all images of $gD$ for all $g\in \Gamma$. This proves that Cayley graph $C'$ has exactly the same set of vertices as the Cayley graph of $D$ and is  connected.  It does not need to be connected in exactly the same way as $C$ and can have less edges joining vertices of $C'$  or few extra edges because $D'$ has less exact faces then $D$ and can few extra faces if $D$ is a proper subset of $D'$. 

A schematic picture is given on the Figure 2 for the 2-dimensional case. The exact Dirichlet domain $D$ is shown in thick solid lines and represents a tiling of a real plane by hexagons, $D'$ is a parallelogram which contains the hexagon $D$ and in addition to $D$ has two shaded triangle. $C'$ is marked with dotted lines. $\square$

\begin{figure}
  \includegraphics[scale=0.3]{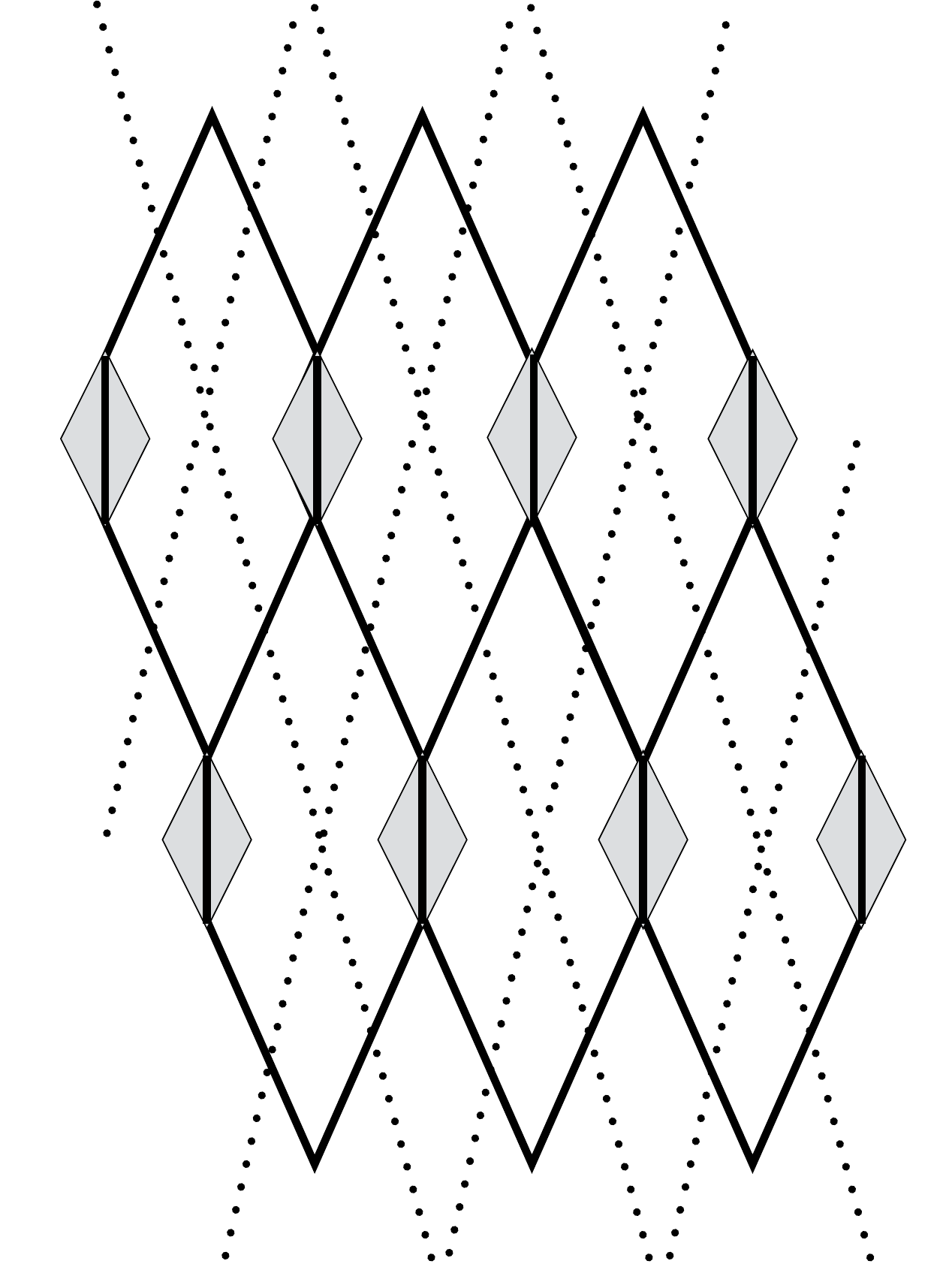}\\
  \caption{Tiling of a real plane by hexagons and its covering by parallelograms.}
\end{figure}

The lemma above tells that the Cayley graph $C'$ is connected in $\mathbb H^3$ and the union of $gD'$ for all $g\in \Gamma$ covers the whole hyperbolic space $\mathbb H^3$ although we replaced the exact Dirichlet domain $D$ with an approximate Dirichlet domain $D'$. The Theorem 1 answers a question if the Cayley graph $C'$ remains connected inside the ball $B(x,R)$.

\medskip

{\bf Proof of Theorem 1:} By construction $D$ and $D'$  share some common faces. Their intersections we call {\it exact walls} and their complements we call {\it hidden walls} in $D$ and {\it extra walls} in $D'$. We want to show that the existence of extra walls does not affect a tiling of the ball $B(x,R)$ by $D'$ if $\Delta V$ is small enough and the spine radius $r'$ is close to minimal. It means that the whole ball $B$ can be tiled by translations of $D'$ if we add neighbors along the existing faces of $D'$ without leaving the ball $B$.  An obstacle we want to avoid here is when a union of domains with basepoints inside $B$ will be disconnected from the $D'_x$ inside the ball $B$ by extra walls of translations of $D'$ although they will remain connected in the universal cover. 

First, we show that an upper bound on the volume difference and an upper bound on the hidden walls area are related. We call it {\it extra area}. If the $\Delta V$ is small enough then the the extra area cannot be very big and cover most of the area of $D$. Assume that $\Delta V$ is small enough but extra area is very big and close to the area of $D$. Then it requires most of faces of $D$ to be almost parallel and very close to the faces of $D'$. This will contradict the fact that images $f_i(x)$ and $f'_i(x)$ of a basepoint $x$ for faces $f_i$ of $D$ and $f'_i$ of $D'$ will belong to different translations of $D$ and will be at distances more than $2\rho$, where $\rho$ is an injectivity radius of $D$. The same argument was used in the proof of Lemma 1.

The first particular case to consider is when there is a single domain $D'_y$ disconnected from other domains in $B$ which intersects the boundary $\partial B$ but the basepoint lies inside $B$ (Figure 3). This $D'_y$ intersects with $B$ along extra walls only. The lower bound on the extra area of $D'$ can be obtained from an area of the inscribed ball in $D'$ intersecting with $B$ divided by a number of neighboring Dirichlet domains $D$ of $D'$. 

\begin{figure}
  \includegraphics[scale=0.45]{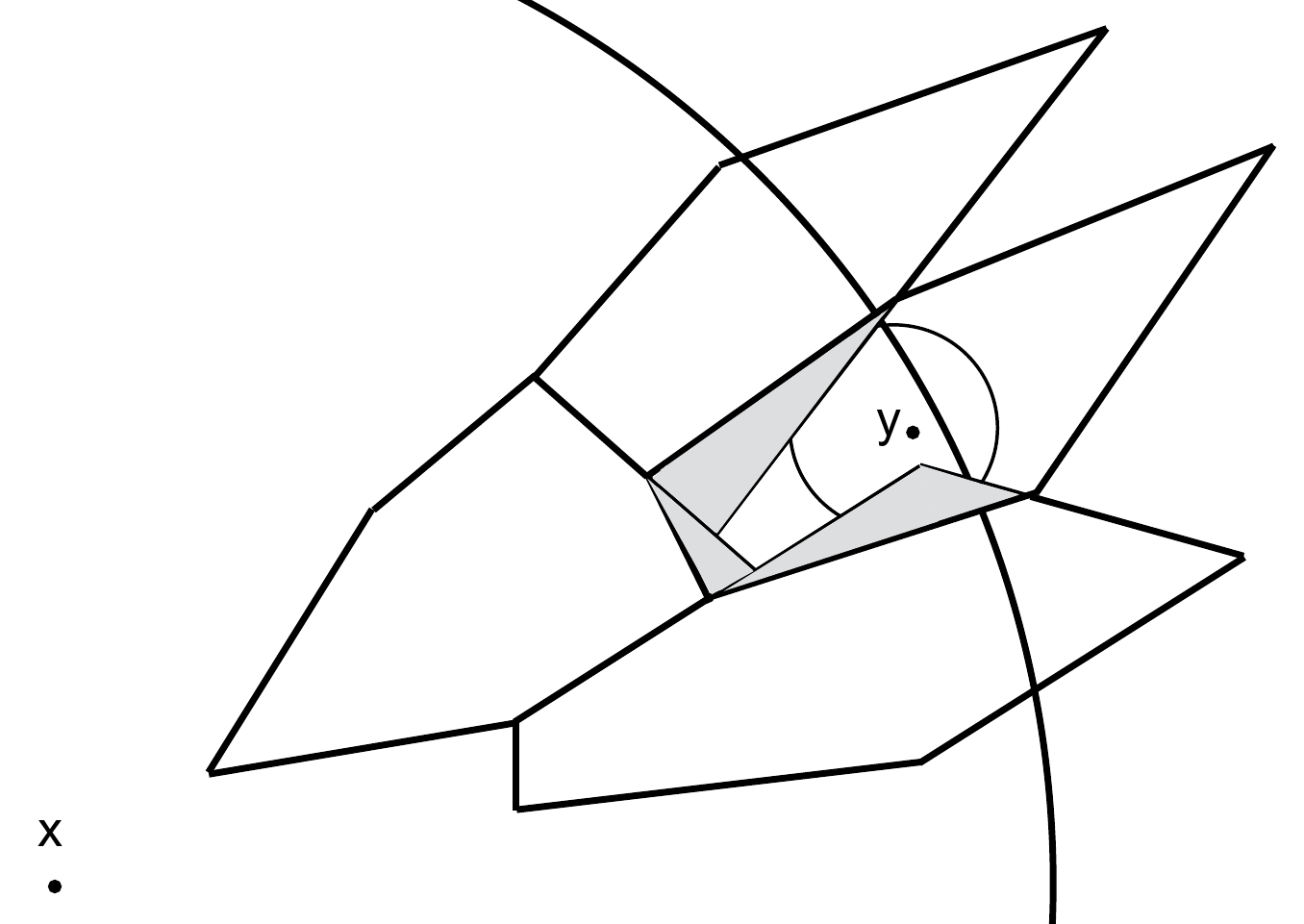}\\
  \caption{Dirichlet domain $D'_y$ is isolated from $D'_x$.}
\end{figure}

This case may not give a lowest bound on the minimal extra area. The extra area will be smaller if it will be intersected by a maximal number of translations of domains $D'$. This can happen in the following case: A single domain $D'_y$ (or a collection of domains $D'$) is surrounded  by several copies of other domains $D'$ and is disconnected in $B$ from $D'_x$ by their extra walls (Figure 4). 

\begin{figure}
  \includegraphics[scale=0.5]{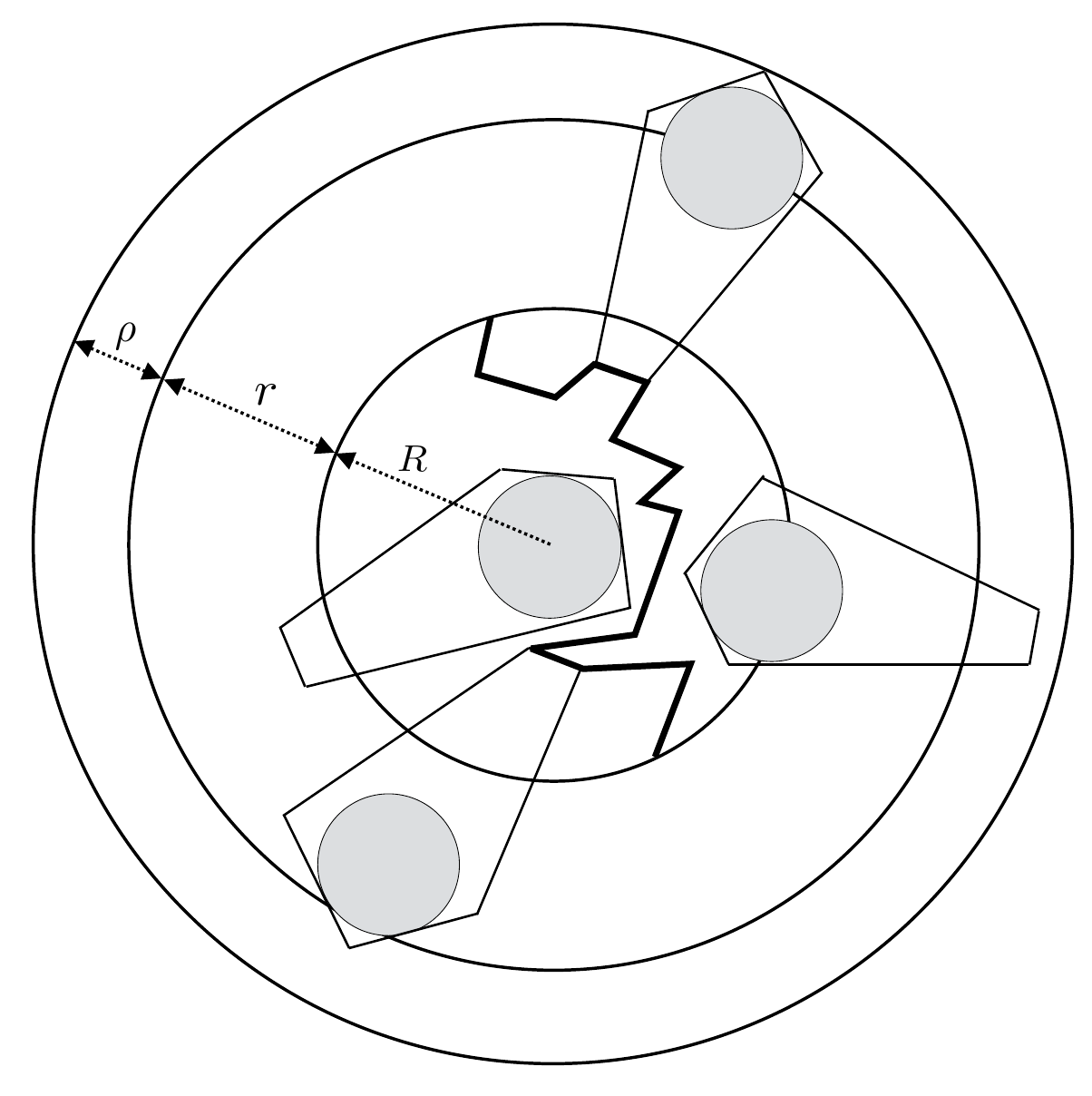}\\
  \caption{Several Dirichlet domains $D'$ are isolated from $D'_x$ in a complicated way. The thick line correspond to a collection of hidden walls.}
\end{figure}

A lower bound of an extra area $``EA"$ is a ratio of a lower bound on an area of the hidden walls $``AHW"$ by an upper bound on a maximum number of Dirichlet domains intersecting these hidden walls $``NDD"$:
$$
\mbox{Extra Area} = \frac{\mbox{Area of Hidden Walls}}{\mbox{Number of intersecting Dirichlet domains}}
$$

The area of a union of these hidden walls $AHW$ is bounded below by the area of the inscribed ball in $D'_y$ intersecting $B$. This area is increasing if a radius of $B$ is increasing. 

It is more difficult to get an upper bound on $NDD$. The essential argument here is the use of a small spine radius. If a Dirichlet domain $D$ has a very long and thin part then an undesired situation can happen when a huge number of translations of $D$ will intersect the hidden walls. In this case a spine radius $r$ is very big. A rough upper bound on the $NDD$ can be obtained by counting a number of inscribed balls in the ball of radius $R+r+\rho$, where $R$ is the radius of the tiling ball $B(x,R)$, $r$ is the spine radius and $\rho$ is the injectivity radius of $D$. The exact value of $r$ cannot be computed from $D'$ but one can get an upper bound of it. It is bounded above by the maximal distance from the basepoint $x$ to the vertices of $D'_x$. A more general case to consider is when there is a union of several Dirichlet domains which are disconnected from $D'_x$ but they all can be reduced to the above case. 

SnapPea's algorithm minimizes the spine radius. It may not construct exactly the minimal spine radius but a one very close to minimal. Therefore, in practice its extra area is much bigger then the one we describe here. An alternative way to explain why in practice an upper bound on NDD cannot be arbitrary large would be by using a maximal injectivity radius $r'$, when an inscribed ball takes most of the volume of $D'$.

{\bf End of proof Theorem 1.}

\section{Computational problem in geometric group theory}

In many situations one needs to compare if two words of a fundamental group $\Gamma\cong\pi_1(M)$ are the same or not. In other words, if two matrices are the same and correspond to the same geodesic in $M$. For example, in the construction of the length spectrum one needs to operate with a big set of matrices which grows exponentially fast. Therefore, all possible duplicates should be removed as earlier as they appear. We present a way  how to solve this word problem rigorously by using approximate data of a Dirichlet domain, namely its injectivity radius.  
 
\medskip
 
 {\bf Lemma 3:} Two elements  $g, g'\in\Gamma$ represent the same element of $\Gamma$ if real parts of their traces are the same and a distance between images $g(x),\,g'(x)$ of the basepoint $x$ is less than twice an injectivity radius of the Dirichlet domain $D_x(\Gamma)$. 

\medskip

{\bf Proof:} Every element $g\in\Gamma$ has a word in terms of generators and their inverses and a matrix that corresponds multiplications of these elements. There are two other numbers which we assign to every $g$: real length of a geodesic and a distance between basepoints $x$ and $g(x)$.

All elements $g,\,g'\in\Gamma$ are different if lengths of their geodesics are different. It is the same as saying that their traces are different. Then we check the difference between the basepoint distances, namely the $\dist(x,g(x))$ and the $\dist(x,g'(x))$. If this difference is more than twice the injectivity radius of $D$ then $g(x)$ and $g'(x)$ can belong to different copies of $D$ and elements  $g, g'$ can be different. If the difference between the basepoint distances is less than the double of the injectivity radius then we have to check another distance: the distance between  basepoints $g(x)$ and $g'(x)$. If this distance is less than the double injectivity radius then $g(x)$ and $g'(x)$ belong to the same Dirichlet domain and $g$ and $g'$ are the same elements of $\Gamma$. 
$\square$

\end{document}